\documentclass{amsart}
\usepackage{amsmath,amsthm,amssymb,amsfonts}
\usepackage{graphicx} 
\usepackage{epstopdf}
\usepackage[T1]{fontenc}
\usepackage[cp1250]{inputenc}
\usepackage{graphicx}
\usepackage{pictex}
\makeatletter
\def\proof{\@ifstar{P\,r\,o\,o\,f}{P\,r\,o\,o\,f.\ }}
\renewcommand\th@remark{%
  \thm@headfont{\bfseries}%
  \normalfont 
  \thm@preskip\topsep \divide\thm@preskip\tw@
  \thm@postskip\thm@preskip
}
\renewenvironment{equation}{\refstepcounter{equation}$$}{\eqno{(\thesection.\theequation)}$$}
\makeatother
\newcounter{Example}[section]
\newcounter{Th}[section] \newcounter{Pr}[section] \newcounter{Lm}[section]
\newcounter{Remark}[section]
\newcounter{Def}[section]
\newcounter{lcounter}[section]
\newcounter{Corol}[section]
\newenvironment{Th}[1][\relax]
    {\medspace\refstepcounter{Th}T\,h\,e\,o\,r\,e\,m \arabic{section}.\theTh.\ \it}
    {\rm\medspace}
\newenvironment{Th.}[1][\relax]
    {\medspace\refstepcounter{Th}T\,h\,e\,o\,r\,e\,m \arabic{section}.\theTh.\ \it}
    {\rm\medspace}
\newenvironment{Pr.}[1][\relax]
    {\medspace\refstepcounter{Pr}P\,r\,o\,p\,o\,s\,i\,t\,i\,o\,n \arabic{section}.\thePr.\ \it}
    {\rm\medspace}
\newenvironment{Lm}[1][\relax]
    {\medspace\refstepcounter{Lm}L\,e\,m\,m\,a \arabic{section}.\theLm.\ \it}
    {\rm\medspace}

\newenvironment{Def}[1][\relax]
    {\medspace\refstepcounter{Def}D\,e\,f\,i\,n\,i\,t\,i\,o\,n \arabic{section}.\theDef.\rm\ }
    {\medspace}
\newenvironment{Def.}[1][\relax]
    {\medspace\refstepcounter{Def}D\,e\,f\,i\,n\,i\,t\,i\,o\,n \arabic{section}.\theDef.\rm\ }
    {\medspace}

\def\au#1{\emph{#1}}
\def\tit#1{{#1}}
\def\R{{\mathbb R}}  
\def\N{{\mathbb N}}  
\def\H{\mathcal{H}} 
\def\cl{\mathop{\rm cl\,}}

\def\co{\mbox{\rm co}\,}

\def\aff {\mbox{\rm aff}\,}

\def\d {\partial\,}
\def\ep{\varepsilon} 
\def\Int {\mbox{\rm int\,}} 

\def\diam {\mbox{\rm diam}\,}

\def\D {\mbox{\rm strco}}
\begin{document}

\title[Uniformly convex subsets with modulus of convexity of the second order]{Uniformly
convex subsets of the
Hilbert space with modulus of convexity of the second order }

\author{Maxim V. Balashov and Du\v{s}an Repov\v{s}}


\address{Department of Higher Mathematics, Moscow Institute of Physics and Technology, Institutski str. 9,
Dolgoprudny, Moscow region, Russia 141700. balashov@mail.mipt.ru}
\address{Faculty of Mathematics and Physics, and Faculty of Education, University of Ljubljana, Jadranska 19, Ljubljana, Slovenia 1000.
dusan.repovs@guest.arnes.si}



\begin{abstract}
We prove that in the Hilbert space every uniformly convex set with
modulus of convexity of the second order at zero is an
intersection of closed balls of fixed radius. We also 
obtain an
estimate of this radius.
\end{abstract}

\date{\today}

\maketitle

\section{Introduction}

\def\i {\mbox{\rm int}\,}

We begin by some definitions for a Banach space $E$ over $\R$.
Let $B_{r}(a)=\{ x\in E\ |\ \| x-a\| \le r\}$. Let $\cl A$ denote
the \it closure \rm and $\Int A$ the \it interior \rm of a subset
$A\subset E$. The \it diameter \rm of a subset $A\subset E$ is
defined as $\diam A = \sup\limits_{x,y\in A} \| x-y\|$.
The \it Minkowski sum \rm of two sets $A,B\subset E$ is the set
$$
A+B = \{ a+b\ |\ a\in A,\ b\in B\}.
$$


We denote the \it convex hull \rm of a set $A\subset E$ by $\co
A$. The \it supporting function \rm of a subset $A\subset E$ is
defined as follows
\begin{equation}\label{**1}
s(p,A) = \sup\limits_{x\in A}(p,x),\qquad\forall p\in E^{*}.
\end{equation}
The supporting function of any set $A$ is always lower
semicontinuous, positively homogeneous and convex. If a set $A$ is
bounded then the supporting function is Lipschitz continuous
\cite{Aubin,Polovinkin+Balashov}. If a subset $A\subset E$ of a
reflexive Banach space $E$ is closed convex and bounded then for
any vector $p\in E^{*}$ the set $A(p)=\{ x\in A\ |\ (p,x)\ge
s(p,A)\}$ is the subdifferential of the supporting function
$s(\cdot,A)$ at the point $p$. In this case the set $A(p)$ is
nonempty, weakly compact and convex (cf. \cite{Aubin,
Polovinkin+Balashov}), $A(0)=A$.

We denote the inner product in the Hilbert space $\H$ by
$(\cdot,\cdot)$.

\begin{Def}\label{modulus} (\cite{Polyak}).
Let $E$ be a Banach space and let a subset $A\subset E$ be convex
and closed. \it The modulus of convexity \rm $\delta_{A}:\
[0,\diam A)\to [0,+\infty) $ is the function defined by
$$
\delta_{A}(\ep) = \sup\left\{ \delta\ge 0\ \left|\
B_{\delta}\left( \frac{x_{1}+x_{2}}{2}\right)\right.\subset A,\
\forall x_{1},x_{2}\in A:\ \| x_{1}-x_{2}\|=\ep \right\}.
$$\rm
\end{Def}

\begin{Def}\label{RM} (\cite{Polyak}).
Let $E$ be a Banach space and let a subset $A\subset E$ be convex
and closed, $A\ne E$. If the modulus of convexity $\delta_{A}(\ep
)$ is strictly positive for all $\ep\in (0,\diam A)$, then we call
the set $A$ \it uniformly convex \rm (\it with the modulus \rm
$\delta_{A}(\cdot)$).\rm
\end{Def}

We proved in \cite{Balashov+Repovs2} that every uniformly convex
set is bounded and if the Banach space  $E$ contains a
nonsingleton uniformly convex set then it admits a uniformly
convex equivalent norm. We also proved that the function $\ep\to
\delta_{A}(\ep)/\ep$ is increasing (see also \cite[Lemma
1.e.8]{Lindestrauss+tzafriri}), and for any uniformly convex set
$A$ there exists a constant $C>0$ such that $\delta_{A}(\ep)\le
C\ep^{2}$.

\begin{Def}\label{UCS} (\cite{Lindestrauss+tzafriri}).
Let $E$ be a Banach space. We call the space $E$ \it uniformly
convex \rm with the modulus $\delta_{E} (\ep)$, $\ep\in [0,2)$,
if the closed unit ball in $E$ is uniformly convex set with the
modulus $\delta_{E}$.
\end{Def}

In a Banach space $E$ consider a set
$$
A = \bigcap\limits_{x\in X}B_{R}(x)\ne\emptyset,
$$
where $X\subset E$ is an arbitrary subset. Such sets have been
considered by several authors (see \cite{Balashov+Polovinkin,
Moreno, Olech+Fr, Polovinkin+Balashov} for details), they are
called \it $R$-convex, \rm or \it strongly convex of radius
$R$\rm. In particular, strongly convex sets of radius $R$ are
closely related to the classical notions of diametrically maximal
sets and constant width sets, see
\cite{Baronti,Bavaud,Martini,Moreno0,Polovinkin,Polovinkin+Balashov}.
It is obvious that if the space $E$ is uniformly convex with the
modulus $\delta_{E}$ then any strongly convex of radius $R$ set
$A$ is uniformly convex with the modulus
$$
\delta_{A}(\ep) \ge R\delta_{E}\left( \frac{\ep}{R}\right),\qquad
\forall \ep\in [0,\diam A).
$$

We want to consider the converse question. Suppose that in a
Banach space $E$ a subset $A\subset E$ is uniformly convex with
the modulus $\delta_{A}$. What can we say about geometric
properties of the set $A$? In particular, is the set $A$ an
intersection of balls of fixed radius?

\section{The Main result}

We give an affirmative answer in the Hilbert space $\H$. Our main
result is given in the following theorem.

\begin{Th}\label{Main} Let $\H$ be the Hilbert space.
Suppose that a nonempty closed convex subset $A\subset\H$ is
uniformly convex with the modulus of convexity of the second order
at zero: there exists $C>0$ such that
$$
\delta_{A}(\ep) = C\ep^{2}+o(\ep^{2}),\quad \ep\to+0.
$$
Then there exists a subset $X\subset\H$ such that
$$
A=\bigcap\limits_{x\in X}\left(x+\frac{1}{8C}B_{1}(0)\right),
$$
and $\frac{1}{8C}$ is sharp in the sense that for any
$r<\frac{1}{8C}$ and any subset $Y\subset\H$,
$$
A\ne \bigcap\limits_{x\in Y}B_{r}(x).
$$
\end{Th}

\section{Preliminary lemmas}

The key idea of the proof of Theorem 2.\ref{Main} is to use the
definition of the \it generating set. \rm

\begin{Def} (\cite{Balashov+Polovinkin})\label{GS}
Let $E$ be a Banach space. A closed convex bounded subset
$M\subset E$ is called a \it generating set, \rm if for any
nonempty subset $A\subset E$ such that
$$
A=\bigcap\limits_{x\in X}(M+x),
$$
where $X\subset E$, there exists another closed convex subset
$B\subset E$ with $A+B=M$. The set $A$ above is said to be \it
$M$-strongly convex.\rm
\end{Def}

\begin{Pr.} \rm (\cite{Balashov+Polovinkin}, \cite[Theorem 4.1.3]{Polovinkin+Balashov}) \label{SP}
Let $M$ be a closed convex bounded set in a reflexive Banach
space $E$. Then $M\subset E$ is a generating set if and only if
for any nonempty set $A=\bigcap\limits_{x\in X}(M+x)$ and any
unit vector $p\in E^{*}$ and any $x_{p}^{A}\in A(p)$ there exists
a point $x_{p}^{M}\in M(p)$ with
\begin{equation}\label{supprinc}
A\subset M+x_{p}^{A}-x_{p}^{M}.
\end{equation}
\end{Pr.}

The inclusion (3.\ref{supprinc}) is a special supporting
principle. Indeed, each closed convex set coincides with the
intersection of supporting half-spaces (\cite{Aubin,
Polovinkin+Balashov}). Proposition 3.\ref{SP} says that if a set
$M$ is generating then each $M$-strongly convex set $A$ coincides
with the intersection of \it supporting shifts \rm of the set $M$.

\begin{Pr.} \rm (\cite{Balashov+Polovinkin}, \cite[Theorem 4.2.7]{Polovinkin+Balashov}) \label{GenerH}
A closed ball in the Hilbert space is a generating set.
\end{Pr.}

We wish to point out that for $\H=\R^{n}$ the result of
Proposition 3.\ref{GenerH} was in fact
proved in \cite{Olech+Fr}.

By Propositions 3.\ref{SP} and 3.\ref{GenerH} we obtain that for
any strongly convex set $A\subset \H$ of radius $R$ and for any
unit vector $p\in \H$ the following
inclusion holds:
\begin{equation}\label{SPinH}
A\subset B_{R}\left(x_{p}^{A}-Rp\right),\qquad \{x_{p}^{A}\} =
A(p).
\end{equation}

\begin{Lm}\label{criteria}
A bounded closed convex set $A\subset \H$ is strongly convex with
the radius $R$ if and only if the function $f(p)=R\| p\|-s(p,A)$
is convex.
\end{Lm}

\proof We conclude from Definition 3.\ref{GS} and Proposition 3.\ref{GenerH}
 that if a closed set $A\subset\H$ is strongly convex
of radius $R>0$ then there exists another convex set $B$ such that
$A+B=B_{R}(0)$. Taking the supporting functions, we get $f(p)=R\|
p\|-s(p,A)=s(p,B)$ which is a convex function.

If the function $f(p)=R\| p\|-s(p,A)$ is convex then, keeping in
mind that $f(p)$ is also continuous and positively homogeneous, we
obtain that $f(p)$ is the supporting function for the set $B=\{
x\in\H\ |\ (p,x)\le f(p),\ \forall p\in\H\}$, i.e. $f(p)=s(p,B)$
\cite[Corollary 1.11.2]{Polovinkin+Balashov}. Hence
$s(p,A)+s(p,B)=s(p,A+B)=R\| p\|$ and by the convexity and
closedness of the set $A$ we have $A+B=B_{R}(0)$. Thus
$A=\bigcap\limits_{b\in B}B_{R}(-b)$.\qed

\begin{Def}\label{strco} (\cite{sca})
For a set $A\subset \H$, $A\subset B_{R}(a)$, a \it strongly
convex hull of radius $R>0$ \rm is defined to be the intersection
of all closed balls of  radius $R$ each of which contains $A$.
We denote the strongly convex hull of radius $R$ of a set $A$ by
$\D_{R}A$.
\end{Def}

Let $\| a-b\|<2R$. Any intersection of the set
$\D_{R}\{a,b\}\subset \H$ by a 2-dimensional plane $L$,
$\{a,b\}\subset L$, represents the planar convex set between two
smaller arcs of the circles of radius $R$ which pass through the
points $a$ and $b$. Also if $0<r<R$ and $\D_{r}A\ne\emptyset$,
then $\D_{R}A\subset \D_{r}A$ (\cite{Balashov+Polovinkin},
\cite[Theorem 4.4.2]{Polovinkin+Balashov}). We define the smaller
arc of a circle of radius $R$, the center $z\in\H$ and the
endpoints $x,y\in\H$ by $D_{R}(z)(x,y)$.

\begin{Lm}\label{local}
Let $R>0$. Let a subset $A\subset \H$ be closed, convex and
bounded. Suppose that
$$
\exists\, \ep_{0}>0\ \forall a,b\in A:\ \| a-b\|\le \ep_{0}\
\Rightarrow\ \D_{R}\{a,b\}\subset A.
$$
Then the set $A$ is strongly convex of radius $R$.
\end{Lm}

\proof The boundary of the set $A$ contains no
nondegenerate 
line
segments.
By the inclusion $\D_{R}\{ a,b\}\subset A$, $\forall a,b\in A$
and $\|a-b \|\le \ep_{0} $, and by the property of strongly convex
hull of two points we obtain that the set  $A$ is uniformly convex
with the modulus
$$
\delta_{A}(t)\ge R-\sqrt{R^{2}-\frac{t^{2}}{4}},\quad\forall t\in
(0,\ep_{0}).
$$
By Corollary 2.2 \cite{Balashov+Repovs2} the function $\d
B_{1}(0)\ni p\to a(p)=\arg\max\limits_{x\in A}(p,x)$ is uniformly
continuous. It is easy to see that $a(p)$ is also uniformly
continuous on each set of the form $\{p\in\H\ |\ \| p\|>r>0\}$.

Fix any pair of linear independent vectors  $p_{1},p_{2}\in\H$
(i.e. $0\notin [p_{1},p_{2}]$). The condition of uniform
continuity on the set  $[p_{1},p_{2}]$ is
\begin{equation}\label{Cantor}
\exists\, \delta_{0}>0\ \forall q_{1},q_{2}\in [p_{1},p_{2}]:\ \|
q_{1}-q_{2}\|<\delta_{0}\qquad \| a(q_{1})-a(q_{2})\|\le \ep_{0}.
\end{equation}
Consider $f(p)=R\| p\|-s(p,A)$, $p\in [p_{1},p_{2}]$. Fix
$q_{1},q_{2}\in [p_{1},p_{2}]$ such that $\|
q_{1}-q_{2}\|<\delta_{0}$. By formula (3.\ref{Cantor}) we obtain
for points $a_{i}=a(q_{i})$, $i=1,2$, that $\|
a_{1}-a_{2}\|\le\ep_{0}$ and using the condition of lemma we have
$\D_{R}\{a_{1},a_{2}\}\subset A$. Using the convexity of the
function $R\| p\| - s(p,\D_{R}\{a_{1},a_{2}\})$ (Lemma
3.\ref{criteria}) we obtain that 

$$
f\left(\frac12 (q_{1}+q_{2})\right)= R\left\| \frac12
(q_{1}+q_{2})\right\|- s\left(\frac12 (q_{1}+q_{2}),A\right)\le
\qquad\qquad\qquad\qquad\qquad\qquad\qquad\qquad\qquad\qquad\qquad\qquad\qquad\qquad\qquad\qquad\qquad\qquad\qquad\qquad\qquad
$$
$$
\qquad\qquad\quad\qquad\le R\left\| \frac12 (q_{1}+q_{2})\right\|
- s\left(\frac12 (q_{1}+
q_{2}),\D_{R}\{a_{1},a_{2}\}\right)\le\qquad\qquad\qquad\qquad\qquad\qquad\qquad\qquad\qquad\qquad\qquad\qquad\qquad\qquad\qquad\qquad\qquad\qquad\qquad\qquad\qquad
$$
$$
\qquad\qquad\quad\qquad\le \frac12 \left( R\|
q_{1}\|-s(q_{1},\D_{R}\{a_{1},a_{2}\})\right)+ \frac12\left( R\|
q_{2}\|-s(q_{2},\D_{R}\{a_{1},a_{2}\})\right)
=\qquad\qquad\qquad\qquad\qquad\qquad\qquad\qquad\qquad\qquad\qquad\qquad\qquad
$$
$$
\qquad\qquad\quad\qquad =\frac12\left( R\|
q_{1}\| - (q_{1},a_{1})\right)+\\
\frac12\left( R\| q_{2}\| -
(q_{2},a_{2})\right)=\qquad\qquad\qquad\qquad\qquad\qquad\qquad\qquad\qquad\qquad\qquad\qquad\qquad\qquad\qquad\qquad\qquad\qquad\qquad\qquad\qquad
$$
$$
\qquad\qquad\quad\qquad=\frac12
f(q_{1})+\frac12f(q_{2}).\qquad\qquad\qquad\qquad\qquad\qquad\qquad\qquad\qquad\qquad\qquad\qquad\qquad\qquad\qquad\qquad\qquad\qquad\qquad\qquad\qquad
$$

Thus for any $q_{1},q_{2}\in [p_{1},p_{2}]$ with $\|
q_{1}-q_{2}\|<\delta_{0}$
$$
f\left( \frac{q_{1}+q_{2}}{2}\right)\le \frac12 f(q_{1})+\frac12
f(q_{2}).
$$

Let us show that the above condition of convexity holds if we
replace $\delta_{0}$ by $2\delta_{0}$, i.e.
$$
\forall q_{1},q_{2}\in [p_{1},p_{2}]:\ \|
q_{1}-q_{2}\|<2\delta_{0}\qquad f\left(
\frac{q_{1}+q_{2}}{2}\right)\le \frac12 f(q_{1})+\frac12 f(q_{2}).
$$
Let $q_{3},q_{4}\in [p_{1},p_{2}]$: $\|
q_{3}-q_{4}\|<2\delta_{0}$, $p_{0}=\frac12 (q_{3}+q_{4})$. Let
$q_{1}=\frac12 (q_{3}+p_{0})$, $q_{2}=\frac12 (q_{4}+p_{0})$, $\|
q_{1}-q_{2}\|<\delta_{0}$. We have $f(p_{0})\le \frac12
(f(q_{1})+f(q_{2}))$, $f(q_{1})\le \frac12 (f(q_{3})+f(p_{0}))$,
$f(q_{2})\le \frac12 (f(q_{4})+f(p_{0}))$. We obtain from the last three
inequalities that $f(p_{0})\le \frac14 f(q_{3})+\frac14
f(q_{4})+\frac12 f(p_{0})$, i.e. $f(p_{0})\le \frac12
f(q_{3})+\frac12 f(q_{4})$.

By induction we obtain that for all $q_{1},q_{2}\in [p_{1},p_{2}]$
$$
f\left(\frac12 (q_{1}+q_{2})\right)\le\frac12 f(q_{1})+\frac12
f(q_{2}).
$$
If $p_{1}$ and $p_{2}$ are parallel then the latter inequality
holds  due to the positive homogeneity of the function $f$.
Finally, by
 continuity of the function $f$ we conclude that $f$ is convex.
Hence, by Lemma 3.\ref{criteria} the set $A$ is strongly convex of
radius $R$.\qed

\begin{Lm}\label{0step}
Let a subset $A\subset\H$ be uniformly convex with the modulus of
convexity $\delta_{A}(\ep)$, $C>0$, and
$\delta_{A}(\ep)=C\ep^{2}+o(\ep^{2})$, $\ep\to+0$. Let $0<K<C$.
Then the set $A$ is strongly convex of radius $\frac{1}{4K}$.
\end{Lm}

\proof By \cite[Theorem 2.1]{Balashov+Repovs2} the set $A$ is
bounded.
 Fix any $K\in (0,C)$. From the asymptotic equality
$\delta_{A}(\ep)\sim C\ep^{2}$, $\ep\to+0$, we obtain that there
exists $\ep_{0}>0$ such  that for all $\ep\in (0,\ep_{0}]$ we
have $\delta_{A}(\ep)>K\ep^{2}$ and
$\delta_{A}(\ep)<\frac{\ep}{2}$. (See Figure 1.)


\begin{figure}[!htb]
\thispagestyle{empty} \font\debeli cmr12 scaled \magstep2
\beginpicture
\setcoordinatesystem units <10 mm, 10 mm> \setplotarea x from
-7.5 to 6, y from -8 to 4.5 \linethickness 0.2mm \putrule from
-4.1312 0 to 4.1312 0 \linethickness 0.3mm \putrule from 0 -8 to
0 4 \put {$m$} [tl] <1mm, 1mm> at 0 4 \put {$w$} [tl] <1mm, -5pt>
at 0 0 \put {$\bullet$} at 0 0 \linethickness 0.15mm \putrule
from -0.3 0 to -0.3 -0.3 \linethickness 0.15mm \putrule from -0.3
-0.3 to 0 -0.3
\circulararc 360 degrees from 2 0 center at 0 0
\put {$z_1$} [bl] <1mm, 1mm> at 0.9683 1.75
\put {$\bullet$} at
0.9683 1.75 \put {$z$} [br] <-1mm, 1mm> at -0.9683 1.75
\put {$\bullet$} at -0.9683 1.75
\setdashpattern <1mm,2mm,1mm,2mm>
\setlinear
\plot 0 0 -0.9683 1.75 /
\setsolid
\setlinear
\plot -5.5771 -0.8 1.5 3.1156 /
\setsolid
\put {$l$} [bl] <1mm, 1mm> at 1.5 3.1156
\setlinear \plot -1.2308 1.6048 -1.0855 1.3423 /
\setsolid
\setlinear
\plot -1.0855 1.3423 -0.8230 1.4875 /
\setsolid
\put {$\bullet$} at -4.1312 0
\put {$a$} [br] <-1mm, 1mm> at -4.1312 0
\setlinear
\plot -4.1312 0 0 -7.467 /
\setsolid
\put {$\bullet$} at 0 -7.4667
\put {$s$} [bl] <1mm, 1mm> at 0
-7.4667
\put {$\bullet$} at 0 2.2857
\put {$x$} [bl] <0.3mm, 2.2mm> at 0 2.2857
\setlinear
\plot -4.3937 -0.1452 -4.2484 -0.4077 /
\setsolid
\setlinear
\plot -4.2484 -0.4077 -3.9859 -0.2625 /
\setsolid
\setlinear
\plot -1.5 3.1156 4.1312 0 /
\setsolid
\put {$l_1$} [br] <-1mm, 1mm> at -1.5 3.1156 \put
{$\bullet$} at 4.1312 0 \put {$b$} [bl] <1mm, 1mm> at 4.1312 0
\circulararc 75 degrees from 4.1312 0 center at 0 -7.4667 \put
{$L\cap \partial B_{\| a-s\|}(s)$} [bl] <1mm, 1mm> at -6.8 -2.2
\put {$L\cap \partial B_{\delta_A(\varepsilon)}(w)$} [tl] <1mm,
1mm> at 1.2 -1.6 \setplotsymbol ({\debeli.}) \circulararc 57.9101
degrees from 4.1312 0 center at 0 -7.4667 \setplotsymbol
({\debeli.}) \circulararc 57.9101 degrees from 0.9683 1.75 center
at 0 0
\endpicture
\medskip

\centerline{Figure 1}
\end{figure}

 Fix an arbitrary pair of points $a,b\in A$ with $\| a-b\|=\ep\le
 \ep_{0}$. Then
 $$
B=\co\left( \{a\}\cup B_{\delta_{A}(\ep)}\left(
\frac{a+b}{2}\right)\cup\{b\}\right)\subset A.
 $$

Consider an arbitrary 2-dimensional affine plane $L$ such that
$\{a,b\}\subset L$. Let $w=\frac12 (a+b)$. Let $l\subset L$ be
 a line such that $a\in l$ and $l$ is a tangent line to the circle
$L\cap\d B_{\delta_{A}(\ep)}(w)$ (at the point $z$). Note that
the segment $[a,z]$ is a part of the boundary $\d B$. Let
$m\subset L$ be a line such that $w\in m$ and $m$ is orthogonal
to the line $\aff\{ a,b\}$. Let the point $s\in m$ be such that
the line $\aff\{ s,a\}$ is orthogonal to the line $l$. Then the
circle $L\cap \d B_{\| a-s\|}(s)$ is tangent to the line $l$ at
the point $a$.

Let the line $l_{1}\subset L$ be symmetric to the line $l$ with
respect to the line $m$. Let the point $z_{1}=l_{1}\cap \left(
L\cap\d B_{\delta_{A}(\ep)}(w)\right)$ be symmetric to the point
$z$ with respect to the line $m$.
Let $R=\| s-a\|$; $\| s-a\|\ge \|
a-w\|=\frac{\ep}{2}>\delta_{A}(\ep)$. Let $x=l\cap m\cap l_{1}$.
The arc $D_{R}(s)(a,b)$ is the homothetic image of the arc
$D_{\delta_{A}(\ep)}(w)(z,z_{1})$ under the homothety with the
center $x$ and the coefficient $k=\frac{\| a-s\|}{\|
z-w\|}=\frac{R}{\delta_{A}(\ep)}$. So we see that
$D_{R}(s)(a,b)\subset L\cap B$.

By the similarity of the triangles $saw$ and $awz$ we have
$\frac{\| z-w\|}{\| a-w\|}=\frac{\| a-w\|}{\| a-s\|}$, or
$$
\frac{2\delta_{A}(\ep)}{\ep}=\frac{\ep}{2R}.
$$
Hence, using inequality $\delta_{A}(\ep)>K\ep^{2}$, we obtain that
$$
R\le \frac{1}{4K}.
$$
By the symmetry of the set $B$ with respect to the line
$\aff\{a,b\}$  and the arbitrary choice of $L$ we have
$$
\D_{\frac{1}{4K}}\{ a,b\}\subset\D_{R}\{ a,b\}\subset \co\left(
\{a\}\cup B_{\delta_{A}(\ep)}\left(
\frac{a+b}{2}\right)\cup\{b\}\right)\subset A.
$$
By Lemma 3.\ref{local} we obtain that the set $A$ is strongly
convex of radius $\frac{1}{4K}$.\qed

\begin{Lm}\label{R1}
Let a subset $A\subset\H$ be uniformly convex with the modulus of
convexity $\delta_{A}(\ep)$, $C>0$, and
$\delta_{A}(\ep)=C\ep^{2}+o(\ep^{2})$, $\ep\to+0$. Let $0<K<C$.
Let the set $A$ be strongly convex of radius $R>\frac{1}{8K}$.
Then the set $A$ is strongly convex of radius
$$
R_{1}=\frac{2R}{8RK+1}.
$$
\end{Lm}

\proof See Figure 2, where $w$, $l$, $m$ etc. are defined as in
proof of Lemma 3.\ref{R1}.

\begin{figure}[!htb]
\thispagestyle{empty}
\font\debeli cmr12 scaled
\magstep2
\beginpicture
\setcoordinatesystem units <13 mm, 13 mm>
\setplotarea x from -6.5 to 5.5, y from -7.5 to 4.5
\linethickness 0.3mm
\putrule from -5.5 0 to 5 0
\linethickness 0.3mm
\putrule from 0 -7.5 to 0 4
\put {$m$} [tl] <1mm, 1mm> at 0 4
\put {$w$} [tr] <-1mm, -5pt> at 0 0
\put {$\bullet$} at 0 0
\linethickness 0.15mm
\putrule from 0.2 0 to 0.2 0.2
\linethickness 0.15mm
\putrule from 0 0.2 to 0.2 0.2
\circulararc 360 degrees from 2 0 center at 0 0
\put {$z_1$} [bl] <1mm, 1mm> at 0.6245 1.9
\put {$\bullet$} at 0.6245 1.9
\put {$\bullet$} at -0.6245 -1.9
\put {$\bullet$} at 0.6245 -1.9
\put {$z$} [br] <-1mm, 1mm> at
-0.6245 1.9 \put {$\bullet$} at -0.6245 1.9
\linethickness 0.15mm
\putrule from -0.6245 -3.5 to -0.6245 3.5
\put {$m_0$} [tr]
<-1mm, 1mm> at -0.6245 3.5
\linethickness 0.15mm
\putrule from
-0.4245 0 to -0.4245 0.2
\linethickness 0.15mm
\putrule from
-0.6245 0.2 to -0.4245 0.2
\linethickness 0.15mm
\putrule from
0.6245 -3.5 to 0.6245 3.5
\put {$m_1$} [tl] <1mm, 1mm> at 0.6245
3.5
\linethickness 0.15mm \putrule from 0.8245 0 to 0.8245 0.2
\linethickness 0.15mm \putrule from 0.6245 0.2 to 0.8245 0.2 \put
{$\bullet$} at -0.6245 -1.9
\put {$\bullet$} at 0.6245 -1.9
\setlinear \plot -0.6245 1.9 2.3008 -7 /
\setsolid \circulararc
40 degrees from -0.6245 1.9 center at 2.3008 -7
\put {$a$} [tr]
<1mm, -2mm> at -3.9256 0
\put {$\bullet$} at -3.9256 0 \setlinear
\plot -4.8150 1 2.3008 -7 / \setsolid
\put {$s$} [bl] <1mm, 1mm>
at 2.3008 -7 \put {$\bullet$} at 2.3008 -7
\put {$c$} [bl] <1mm, 1mm> at 0 -4.4133
\put {$\bullet$} at 0 -4.4133
\setlinear
\plot -5.2747 -1.2 -1.4522 2.2 /
\setsolid
\put {$l$} [br] <-2mm, -0.5mm> at -1.4522 2.2
\setlinear
\plot -4.0585 0.1494 -3.9090 0.2824 /
\setsolid
\setlinear
\plot  -3.9090 0.2824 -3.7761 0.1329 /
\setsolid
\circulararc -83.3047 degrees from -3.9256 0
center at 0 -4.4133
\put {$b$} [bl] <1mm, 1mm> at 3.9256 0
\put {$\bullet$} at 3.9256 0
\put {$x$} [tr] <-1mm, -1mm> at -0.6245
1.4601
\put {$\bullet$} at -0.6245 1.4601
\put {$x_1$} [tl] <1mm, -1mm> at 0.6245 1.4601
\put {$\bullet$} at 0.6245 1.4601
\arrow <2mm> [0.2,0.7] from 3 2 to 1.75 1.25
\setplotsymbol ({\debeli.})
\circulararc -23.4575 degrees from 0.6245 1.9 center at -2.3008 -7
\circulararc 23.4575 degrees from 0.6245 -1.9 center at -2.3008 7
\circulararc -23.4575 degrees from -0.6245 -1.9 center at 2.3008 7
\circulararc 23.4575 degrees from -0.6245 1.9 center at 2.3008 -7
\put {$L\cap\partial B_{R}(s)$} [bl] <1mm, 1mm> at -6 -2.7 \put
{$D_1$} [bl] <1mm, 1mm> at -1.9 1.5
\put {$D_2$} [br] <1mm, 1mm>
at 1.9 1.5 \put {$D_3$}  at -2.3 -0.9
\put {$D_4$}  at 2.3 -0.9
\put {$D_{\varrho}(c)(a,b)$} [bl] <1mm, 1mm> at 3 2 
\endpicture
\medskip
\centerline{Figure 2}
\end{figure}

Fix any $K\in (0,C)$. From the asymptotic equality
$\delta_{A}(\ep)\sim C\ep^{2}$, $\ep\to+0$, we obtain that there
exists such $\ep_{0}\in\left(0, \frac{1}{100K}\right)$
 that for all $\ep\in (0,\ep_{0}]$ we have
 $\delta_{A}(\ep)<\frac{\ep}{2}$ and
 \begin{equation}\label{good_local}
\frac{\delta_{A}(\ep)}{\ep^{2}}\left(
 1-\frac{\delta_{A}(\ep)}{R}\right)>K.
  \end{equation}

Choose any pair of points $a,b\in A$, $\| a-b\|=\ep\le \ep_{0}$.
Then
$$
B = \D_{R}\left\{ \{a\}\cup B_{\delta_{A}(\ep)}(w)\cup\{ b
\}\right\}\subset A,
$$
where $w=\frac{a+b}{2}$.

Let $L$ be any 2-dimensional affine plane, $\{ a,b\}\subset L$.
 Let $m\subset L$ be a line such that $w\in
m$ and $m\perp \aff\{a,b\}$. Consider a circle of radius $R$ with
a center $s\in L$ such that $L\cap B_{R}(s)\supset L\cap
B_{\delta_{A}(\ep)}(w)$, it passes through the point $a$,
tangents to $L\cap\d B_{\delta_{A}(\ep)}(w)$ and define
$$
\left[ L\cap \d B_{R}(s)\right]\cap \left[ L\cap\d
B_{\delta_{A}(\ep)}(w)\right]=\{ z\}.
$$
 Such
circle exists because $R>\ep/2$.

Let $D_{1}=D_{R}(s)(a,z)$, $D_{2}$ is symmetric to the $D_{1}$
with respect to the line $m$, $D_{3}$ and $D_{4}$ are symmetric to
the $D_{1} $ and $D_{2}$ with respect to the line $\aff\{ a,b\}$,
respectively. We have
$$
L\cap B \supset \co\left\{ D_{1}\cup D_{2}\cup D_{3}\cup D_{4}\cup
\left(L\cap B_{\delta_{A}(\ep)}(w)\right) \right\}.
$$

Let $l$ be the tangent line to the circle $L\cap\d B_{R}(s)$ at
the point $a$. Let $\varphi$ be the angle between the lines $l$
and $\aff\{ a,b\}$; $\alpha$ be the angle between the lines
$\aff\{ a,b\}$ and $\aff\{ a,s\}$; $\varphi+\alpha=\pi/2$.
Consider the triangle $aws$: $\| a-w\|=\frac{\ep}{2}$, $\|
a-s\|=R$, $\| w-s\| = R-\delta_{A}(\ep)$. Hence by the Cosine
theorem we get
$$
\sin\varphi = \cos\alpha  = \frac{\| a-s\|^{2}+\| a-w\|^{2}-\|
w-s\|^{2}}{2\| a-s\|\cdot\| a-w\|} =
$$
$$
=\frac{\ep}{4R}+\frac{2\delta_{A}(\ep)}{\ep}-\frac{\delta^{2}_{A}(\ep)}{R\ep}
=\frac{\ep}{4R}+\frac{2\delta_{A}(\ep)}{\ep}\left(
1-\frac{\delta_{A}(\ep)}{R}\right),
$$
and using (3.\ref{good_local}) we obtain that
\begin{equation}\label{sin+cos}
\sin\varphi = \cos\alpha \ge \frac{\ep}{4R} + 2K\ep.
\end{equation}

Let the point $z_{1}=D_{2}\cap\left( L\cap \d
B_{\delta_{A}(\ep)}(w)\right)$ be symmetric to the point $z$ with
respect to the line $m$. Let the lines $ m_{0}\subset L$ and
$m_{1}\subset L$ be parallel to the line $m$, $z\in m_{0}$,
$z_{1}\in m_{1}$. Let $\varrho = \frac{\ep}{2\sin \varphi}$ and
$c=m\cap\aff\{ a,s\}$. Note that $\| a-c\|=\varrho$. Let
$x=D_{\varrho}(c)(a,b)\cap m_{0}$, $x_{1}=
D_{\varrho}(c)(a,b)\cap m_{1}$.

The circle $L\cap \d B_{\varrho}(c)$ touches the line $l$ at the
point $a$. Taking into account that $K>1/8R$, we get by the
formula (3.\ref{sin+cos})
$$
\varrho = \frac{\ep}{2\sin\varphi}\le\frac{1}{\frac{1}{2R}+4K}<
\frac{1}{\frac{1}{2R}+\frac{1}{2R}}=R.
$$
By the last inequality the points $a$ and $s$ are separated by the
line $m$ and the points $a$, $z$ (and $b$, $z_{1}$) are situated
in the same halfplane with respect to the line $m$.  Hence the
points $x$ and $x_{1}$ belong to the disk $L\cap
B_{\delta_{A}(\ep)}(w)$. The arc $D_{\varrho}(c)(a,x)$ lies
between the arc $D_{R}(s)(a,z)$ and the line $\aff\{ a,b\}$, by
the symmetry of the arc $D_{\varrho}(c)(a,b)$ with respect to the
line $m$ we have
$$
D_{\varrho}(c)(a,x)\cup D_{\varrho}(c)(b,x_{1})\subset L\cap B.
$$
 By the inequality $\varrho\ge
\frac{\ep}{2}>\delta_{A}(\ep)$ we have
$$
D_{\varrho}(c)(x,x_{1})\subset L\cap B_{\delta_{A}(\ep)}(w).
$$
Thus $D_{\varrho}(c)(a,b)\subset L\cap B$ and $\varrho\le
\frac{1}{\frac{1}{2R}+4K}=R_{1}$. By the symmetry of the set
$L\cap B$ with respect to the line $\aff\{ a,b\}$ and the
arbitrary choice of $L$ we have
$$
\D_{R_{1}}\{a,b\}\subset \D_{\varrho}\{a,b\}\subset B\subset A.
$$
By Lemma 3.\ref{local} we obtain that the set $A$ is strongly
convex with the radius $R_{1}$.\qed

\begin{Lm}\label{last}
Under the assumptions of Lemma 3.\ref{R1} the set $A$ is strongly
convex of radius $r=\frac{1}{8K}$.
\end{Lm}

\proof Define inductively the sequence $\{
R_{n}\}_{n=1}^{\infty}$ as follows: $R_{1}=R$, and, for $n\ge 1$,
$$
R_{n+1}=\frac{2R_{n}}{8R_{n}K+1}.
$$
It is not difficult to show that $R_{n+1}\le R_{n}$. After having
a look to the function $f(x)=\frac{2x}{8xK+1}$ (which is
increasing for $x\ge 0$), and using $R>\frac{1}{8K}$, it is not
difficult to show that also $R_{n}>\frac{1}{8K}$ for every
$n\in\N$. Hence by the Weierstrass theorem we get $R_{n}\to
r=\frac{2r}{8rK+1}$, that is, $r=\frac{1}{8K}$.

Using Lemma 3.\ref{R1} we know that the set $A$ is strongly convex
of radius $R_{n}$ for every $n\in \N$. By Lemma 3.\ref{local},
the latter assertion means that the functions $f_{n}(p)=R_{n}\|
p\|-s(p,A)$ are convex for all $n\in\N$. Taking the limit of the
sequence of the functions we get that
$$
f(p)=\frac{1}{8K}\| p\|-s(p,A)
$$
is a convex function as well. This shows, again by Lemma
3.\ref{local}, that $A$ is strongly convex of radius
$\frac{1}{8K}$.\qed

\section{Proof of Theorem 2.\ref{Main}}

 Let $K_{n}=C-\frac1n$. By
Lemma 3.\ref{0step} the set $A$ is strongly convex of radius
$\frac{1}{4K_{n}}$. By Lemma 3.\ref{last} the set $A$ is strongly
convex of radius $\frac{1}{8K_{n}}$. By Lemma 3.\ref{criteria}
this is equivalent to the convexity of the function
$\frac{1}{8K_{n}}\| p\| - s(p,A)$ for all natural $n$. Taking the
limit $n\to \infty$, we obtain the convexity of the function
$\frac{1}{8C}\| p\| - s(p,A)$. Hence by Lemma 3.\ref{criteria}
the set $A$ is strongly convex of radius $\frac{1}{8C}$.

Suppose that there exist a number $r\in\left(
0,\frac{1}{8C}\right)$ and a subset $Y\subset \H$ such that
$$
A = \bigcap\limits_{x\in Y}B_{r}(x).
$$
 For the ball
$B_{r}(0)$ in the Hilbert space we have (cf.
\cite{Lindestrauss+tzafriri})
$$
\delta_{B_{r}(0)}(\ep)=r\delta_{\H}\left(\frac{\ep}{r}\right) =
r-\sqrt{r^{2}-\frac{\ep^{2}}{4}}=
\frac{\ep^{2}}{8r}+o(\ep^{2}),\quad \ep\to+0.
$$
Due to the fact that the set $A$ is the intersection of closed
balls of the radius $r$ we conclude
that
$C\ep^{2}+o(\ep^{2})=\delta_{A}(\ep)\ge
r\delta_{\H}\left(\frac{\ep}{r}\right)=
\frac{\ep^{2}}{8r}+o(\ep^{2})$, for all $\ep>0$. Hence $C\ge
\frac{1}{8r}$. This contradicts the inequality
$r<\frac{1}{8C}$.\qed

 \section*{Acknowledgements}

  This research was supported by SRA grants P1-0292-0101,
 J1-2057-0101, and BI-RU/10-11-002. The first author
was supported by  RFBR grant 07-01-00156, ADAP project
"Development of scientific potential of higher school" 2.1.1/500
and and project of FAP "Kadry" 1.2.1 grant P938 and grant
16.740.11.0128.
We thank the referee for several comments and suggestions.

\end{document}